\documentclass{elsart}
\usepackage{amssymb}
\usepackage{amsfonts}
\usepackage{amsmath}
\usepackage{epsf}
\usepackage{epsfig}
\usepackage{graphicx}

\input{tcilatex}
\begin{document}

\begin{frontmatter}%

\title{Controlling infectious diseases: the decisive phase effect on a seasonal vaccination strategy}%

\author{J. Duarte$^1$$^,$$^*$, C. Janu\'ario$^2$, N. Martins$^3$, J. M. Seoane$^4$, M. A. F. Sanjuan$^4$}%

\address{$^1$ISEL-Engineering Superior Institute of Lisbon, Department of Mathematics,
\\Rua Conselheiro Em\'idio Navarro 1, 1950-007 Lisboa and Center for Mathematical Analysis, 
Geometry and Dynamical Systems, Instituto Superior T\'ecnico, Universidade de Lisboa, Portugal
\\$^2$Center for Research and Development in Mathematics and Applications (CIDMA) 
Department of Mathematics, University of Aveiro, 3810-193 Aveiro, Portugal
\\$^3$Department of Mathematics and Center for Mathematical Analysis, Geometry and Dynamical Systems, 
Instituto Superior T\'ecnico, Universidade de Lisboa, Av. Rovisco Pais 1, 1049-001 Lisboa, Portugal
\\$^4$Nonlinear Dynamics and Chaos Group, Departamento de F\'{\i}sica, 
Universidad Rey Juan Carlos, Tulip\'{a}n s/n 28933 M\'{o}stoles, Madrid, Spain
\\$^*$E-mail: jorge.duarte@isel.pt (corresponding author)}%

\begin{abstract}%
The study of epidemiological systems has generated deep interest in
exploring the dynamical complexity of common infectious diseases driven by
seasonally varying contact rates. Mathematical modeling and field
observations have shown that, under seasonal variation, the incidence rates
of some endemic infectious diseases fluctuate dramatically and the dynamics
is often characterized by chaotic oscillations in the absence of specific
vaccination programs. In fact, the existence of chaotic behavior has been
precisely stated in the literature as a noticeable feature in the dynamics
of the classical Susceptible-Infected-Recovered (SIR) seasonally forced
epidemic model. However, in the context of epidemiology, \textit{chaos} is
often regarded as an undesirable phenomenon associated with the
unpredictability of infectious diseases. As a consequence, the problem of
converting chaotic motions into regular motions becomes particularly
relevant. In this article, we consider the phase control technique applied
to the seasonally forced SIR epidemic model to suppress chaos.
Interestingly, this method of controlling chaos takes on a clear meaning as
a weak perturbation on a seasonal vaccination strategy. Numerical
simulations show that the phase difference between the two periodic forces -
contact rate and vaccination - plays a very important role in controlling
chaos.%
\end{abstract}%


\end{frontmatter}%

\section{Introduction}

The outbreaks and spread of infectious diseases are a threat to public
health and a source of serious problems to the economic and social
development of our society. The research in epidemic dynamics is, without
doubt, critical to any attempt to prevent or minimize the transmission of
diseases. In particular, we have all witnessed a remarkable intensification
of this work in the present days given the quantity of infected people
worldwide due to the dissemination of a new virus, referred to as COVID-19
(Coronavirus disease 2019), an infectious disease emerged from China in
November 2019 \cite{Libotte2019}.

The use of mathematical models has contributed greatly to our better
understanding of the underlying mechanisms that influence the spread of
diseases and has suggested control strategies, which may have significant
management implications \cite{Zhang2013}. More precisely, these models are
truly significant in different fields, such as policy making, risk
assessment, emergency planning and definition of health-economic
control-programs. The mechanism of transmission of infections is now known
for most diseases. Generally, diseases transmitted by viral agents (such as
measles, influenza, rubella and chickenpox) confer immunity against
reinfection, while diseases transmitted by bacteria (such as meningitis and
tuberculosis) confer no immunity against reinfection. Other diseases, such
as malaria, are transmitted not directly from human to human but by vectors
(usually insects). The goal of mathematical epidemiology has been the
development of mathematical models for the spread of disease as well as
tools for their analysis.

The first mathematical model describing an infectious disease was proposed
by Bernoulli, in 1760, to study the spread of smallpox, which was prevalent
at the time \cite{Bernoulli1760}. An early recognition of the importance of
mathematical modeling came in 1911 when Ronald Ross won the Nobel Prize in
Physiology or Medicine, for his demonstration of the dynamics of the
transmission of malaria between mosquitoes and humans (using models
formulated with differential equations) \cite{Ross1911}. Recent studies in
epidemiology are particularly focused on modeling the previously mentioned
COVID-19 and high impact recurrent epidemics, best exemplified by childhood
infectious diseases such as measles, chickenpox, mumps, whooping cough and
rubella, but also including hepatitis, different types of influenza (see 
\cite{Stone2008} and references therein).

One of the questions that first attracts the attention of scientists
interested in the study of the spread of communicable diseases is how severe
will an epidemic be. This question may be interpreted in a variety of ways
and studied with the aid of models. To formulate dynamic models predicting
the behavior of outbreaks and the transmission of infectious diseases,
compartmental models are usually considered, with the population in a given
region subdivided into several distinct groups or compartments, depending
upon their experience with respect to disease, whose sizes change with time.
These compartments, which are mutually exclusive categories based on
infection status, were initially proposed by William O. Kermack and Anderson
G. McKendrick in 1927 with their mathematical epidemic nonlinear system of
differential equations, called the Susceptible-Infected-Recovered (SIR)
model \cite{Mckendrick1927}. In the context of this prominent model, of
great historical importance in research of epidemics, the population being
studied is labeled into three classes: the susceptible compartment refers to
individuals who have never come into contact with the disease at time $t$,
but they could catch it, i.e., they are vulnerable to exposure with
infectious people. Infected individuals, assumed infectious and capable of
spreading the disease to those in the susceptible category, and remain in
the infectious compartment until their recovery. The recovered class refers
to individuals who have been infected and then recovered from the disease.
The recovered individuals are immune for life and are not able to transmit
the infection to others. They are essentially removed from the population
and play no further role in the dynamics. These models can include, among
other aspects, time-dependent parameters to represent the effects of
seasonality and human demographics by adopting birth and death rates. In
practice, given a model, its parameters must be determined to represent a
particular epidemic context. In formulating models as systems of
differential equations, we are assuming that the epidemic process is
deterministic, that is, the behavior of a population is determined
completely by its history and by the rules which describe the model.
Particularly, in the present study, we will consider a SIR model
representing a class of seasonally forced epidemic models with vital
dynamics (birth and death rates) and constant population.

The article is organized as follows. Immediately after the present Section
I, where an introduction to our study is provided, Section II gives
essential preliminaries about the model, as well as the control procedure
used. Section III is devoted to a diagnosis of the existence of chaos in the
system. We exhibit different dynamical interactions between the densities of
the susceptible and the infected individuals varying the degree (or
amplitude) of seasonality, $\varepsilon $. Revealing insights about the
chaotic dynamics of the SIR model, incorporating seasonal fluctuation, are
gained through the computation of bifurcation diagrams and the largest
Lyapunov exponent. In Section IV, the phase control technique is applied in
order to suppress chaos. A clear biological meaning is assigned to this
method of controlling chaos as a weak/small perturbation on a seasonal
vaccination strategy. Finally, our last considerations are devoted to
significative conclusions.

\section{Description of the model}

In the literature, and as far as the description of the model is concerned,
two related sets of time-dependent variables have been considered in the
epidemic modeling process. The first set of dependent variables counts the
number of people in each of the groups, each as a function of time.
Considering that the population size ($N$) is defined as $N=\widetilde{S}+%
\widetilde{I}+\widetilde{R}$, the dynamical behavior of a homogeneously
mixing population is described using the compartmental forced $\widetilde{S}%
\widetilde{I}\widetilde{R}$ model \cite{Dietz1976}, with vital dynamics
(birth and death rates), given by

\begin{equation}
\left \{ 
\begin{array}{l}
\frac{d\widetilde{S}}{dt}=\sigma -\mu \widetilde{S}-\beta (t)\frac{%
\widetilde{S}~\widetilde{I}}{N} \\ 
\frac{d\widetilde{I}}{dt}=\beta (t)\frac{\widetilde{S}~\widetilde{I}}{N}%
-(\gamma +\mu )\widetilde{I} \\ 
\frac{d\widetilde{R}}{dt}=\gamma \widetilde{I}-\mu \widetilde{R}%
\end{array}%
\right. .  \label{SIR_Model0}
\end{equation}%
The second set of dependent variables represent a fraction of the total
population in each of the three categories. So, being $N$ the previous total
population, we have

\begin{eqnarray}
S &=&\frac{\widetilde{S}(t)}{N}\text{, the susceptible fraction of the
population,}  \label{SIRTilda} \\
I &=&\frac{\widetilde{I}(t)}{N}\text{, the infected fraction of the
population,}  \notag \\
R &=&\frac{\widetilde{R}(t)}{N}\text{, the recovered fraction of the
population,}  \notag
\end{eqnarray}%
verifying $S+I+I=1$. The two sets of dependent variables are proportional to
each other, giving the same information about the progress of the epidemic.
The independent variable is time $t$, measured in a specific unit according
to the epidemic context (such as days, years, etc.). Although it may seem
more natural to work with the populations counts, some of our calculations
will be simpler if we use the fractions instead. As a consequence, we are
going to consider the $SIR$ model with vital dynamics and constant
population $(\forall t,$ $\frac{dS}{dt}+\frac{dI}{dt}+\frac{dR}{dt}=0)$,
such that $S+I+R=1,$%
\begin{equation}
\left \{ 
\begin{array}{l}
\frac{dS}{dt}=\sigma -\mu S-\beta (t)\frac{SI}{N} \\ 
\frac{dI}{dt}=\beta (t)\frac{SI}{N}-(\gamma +\mu )I \\ 
\frac{dR}{dt}=\gamma I-\mu R%
\end{array}%
\right. ,  \label{SIR_Model1}
\end{equation}%
with time $t$\ scaled in unit of years. The parameters of the model (\ref%
{SIR_Model1}) and respective meaning are presented in the following table
(for more details, please see \cite{Zhang2013}, \cite{Schwartz1983}, \cite%
{Kamo2002} and references therein).

\begin{center}
Table 1. Description of the parameters of model (\ref{SIR_Model1})

$%
\begin{tabular}{|l|l|}
\hline
Parameter (year)$^{-1}$ & Description \\ \hline
$\sigma $ & Birth rate \\ \hline
$\mu $ & Natural death rate \\ \hline
$\gamma $ & R$\text{ecovery rate}$ \\ \hline
$\beta (t)$ & Contact or transmission rate \\ \hline
\end{tabular}%
$
\end{center}

We assume that the newborns are susceptible and that the birth and death
rates are balanced, $\sigma =\mu $. As we will exemplify with particular
values of the parameters, the average time to recover from infection can be
derived from $\gamma $\ and given by $1/\gamma $. The contact or
transmission rate of the infection (the coefficient of infectivity) $\beta
(t)$ represents the number of contacts with other individuals per infective
per unit of time. More intensive research on seasonally-forced epidemic
models did not begin until the 1970s (\cite{London1973}, \cite{Schenzle1984}%
). With the purpose of expressing seasonality, a commonly used scheme for
the contact rate $\beta (t)$ takes the form%
\begin{equation}
\beta (t)=\beta _{0}\left( 1+\varepsilon \varphi (t)\right) ,
\label{GeneralBeta}
\end{equation}%
where $\beta _{0}$ gives the mean contact rate; parameter $\varepsilon $,
with $0\leq \varepsilon \leq 1$, represents the strength of the seasonal
forcing (measuring the degree of seasonality); $\varphi $ is a $T$ -
periodic function of zero mean and $t$ is scaled in units of years. Modelers
often incorporate seasonality by making the contact rate $\beta (t)$ a
sinusoidal function of time $\beta _{\cos }(t)=\beta _{0}\left(
1+\varepsilon \cos \left( 2\pi t\right) \right) $. With this procedure, a
seasonally-varying transmission rate yields oscillations at periods that are
integer multiples of the period of forcing.

\begin{figure}[th]
\center{\includegraphics[width=0.45\textwidth]{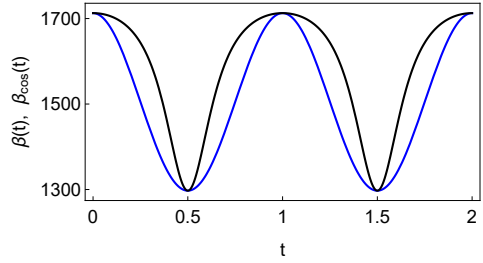}}
\caption{Comparison between the periodic Kot-type function $\protect \beta %
(t) $ (Black) and the usual periodic function $\protect \beta _{\cos }(t)$
(Blue) ($\protect \varepsilon =0.138$).}
\label{BetaSinusoidal}
\end{figure}

As recently proposed in \cite{Buonomo2017}, and inspired by the developed
arguments, we are going to consider in our epidemiological study yearly
periodic Kot-type functions. Particularly, the previously mentioned function 
$\beta (t)$ is given by

\begin{equation}
\beta (t)=\beta _{0}\left( 1+\varepsilon \left( \frac{\frac{2}{3}+\cos
\left( 2\pi t\right) }{1+\frac{2}{3}\cos \left( 2\pi t\right) }\right)
\right) \text{.}  \label{ParticularBeta}
\end{equation}%
In order to point out the main features of the adopted type of periodic
functions, we compare in Fig.~\ref{BetaSinusoidal} a sinusoidal-type
periodic function with the periodic Kot-type function $\beta (t)$.
Differently from the periodic sinusoidal-type function, the Kot-type
function gives asymmetrical weights to the seasonal regimes, stressing the
relative maxima of the contact rate. This is an eye-catching and noteworthy
feature that has been pointed out as more realistic.

An exposure of a susceptible to an infectious is an encounter in which the
infection is transmitted. In this context, the contact rate, $\beta (t)$, is
the average density of susceptible in a given population contacted per
infectious individuals per unit of time. Therefore, $\beta (t)S(t)$ denotes
the rate of the total density of susceptible infected by one infectious and $%
\beta (t)S(t)I(t)$ represents the rate of infection of the susceptible by
all infectious.

It is well known that seasonal forces, including climatic factors and human
phenomena (such as school schedules), are a primary factor responsible for
the transmission and dynamical behavior of most recurrent infectious
diseases \cite%
{Yorke1973,Clarkson1982,Keeling1997,Earn2000,Stone2007,Diedrichs2014,Axelsen2014}%
. As a consequence, it becomes natural to model these diseases as
periodically forced nonlinear systems (\cite%
{Stone2007,Billings2002,Olsen1990,Aron1984,Keeling2001}). Seasonal forces in
these systems shape the spread of the infectious disease and many studies
have stated that intense seasonality can lead to strong erratic patterns
with the presence of chaotic oscillations (\cite{Herrera2017} and \cite%
{Olsen1990}). The chaotic behavior of the simple SIR model with seasonality
has deep biological consequences. The sensitive dependence on the initial
conditions, in the final epidemic outcome, is one of the main concerns. It
is important to notice that this rich dynamical scenario is a direct
consequence of seasonality, which makes any vaccination strategy difficult
to design. In fact, the inclusion of seasonality increases the complexity of
the models at several levels, making a complete comprehension of its
influence a non-trivial and challenging task. In this context, chaos is
regarded as an unwanted feature and, as a consequence, control schemes
(considered as methods of controlling chaos, suppressing chaos, or taming
chaos) play a prominent/decisive role. In any case, the explicit aim of all
the procedures is to obtain stable periodic orbits from chaotic ones by
applying a small/weak, and carefully chosen, perturbation to the system.

The controlling methods have been traditionally classified into two general
categories: feedback methods and nonfeedback methods, based on their
interaction with the system's dynamics. Feedback methods attempt to suppress
chaos by stabilizing orbits already existing in the chaotic attractor.
Nonfeedback methods have been essentially used to stabilize periodically
driven chaotic dynamical systems. Typically, they make stable solutions to
appear by applying small driving forces directly to some of the parameters
of the system or as an additional forcing. It has been shown in the
literature (for instance, in papers \cite%
{Yang1996,Zambrano2006,Zambrano2008,Joseph2012,Used2012}, and references
therein) that a phase difference $\phi $, between the main driving and the
perturbation, influences profoundly the global dynamics of the system.

As we will see in the following lines, given the nature of the forced SIR
epidemic model, it is specifically a nonfeedback procedure, that makes use
of the property where $\phi$ acts as a control parameter, that stands out to
be particularly useful in controlling chaotic behavior. Acting as a
positive/beneficial biological consequence, this nonfeedback scheme called
phase control of chaos makes a vaccination strategy clearly easier to design.

Following the previous considerations, the particular system we are going to
consider in this article is%
\begin{equation}
\left \{ 
\begin{array}{l}
\frac{dS}{dt}=\sigma -\mu S-\beta (t)SI-v\left( t\right) S \\ 
\frac{dI}{dt}=\beta (t)SI-(\gamma +\mu )I \\ 
\frac{dR}{dt}=\gamma I-\mu R+v\left( t\right) S%
\end{array}%
\right. ,  \label{SystemVaccination}
\end{equation}%
where the first periodic force, the contact rate $\beta (t)$, drives the
system to the chaotic state, while the second one, the vaccination rate $%
v\left( t\right) $, is a weak periodic perturbation sensitively modifying
the system's dynamics. We assume that the vaccination rate $v\left( t\right) 
$ is defined by the Kot-type function

\begin{equation}
v(t)=v_{0}+\alpha \left( \frac{\frac{2}{3}+\cos \left( 2\pi rt+\phi \right) 
}{1+\frac{2}{3}\cos \left( 2\pi rt+\phi \right) }\right) ,
\label{Additive_Perturbation}
\end{equation}%
where $v_{0}$ gives the mean vaccination rate, $\phi $ is the phase
difference between the applied perturbation $v\left( t\right) $ and the
driving force $\beta (t)$, a key parameter for our control scheme. Parameter 
$r$ is the ratio of the frequency of those forces and $\alpha $, with $%
\alpha <<1$, measures the degree (or amplitude) of the seasonality of $%
v\left( t\right) $. Once vaccinated, the individuals are no more susceptible
(fact represented by the term $-v\left( t\right) S$\ in the first equation
of the model) and the vaccinated individuals are added to the recovered
class (fact represented by the term $+v\left( t\right) S$\ in the third
equation of the model).

\begin{figure}[th]
\center{\includegraphics[width=0.45\textwidth]{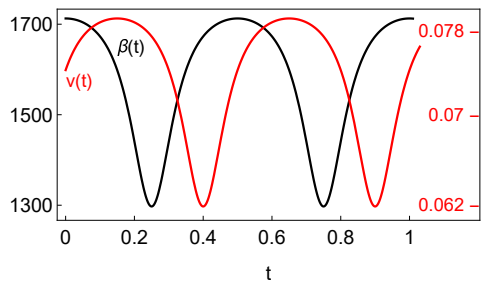}}
\caption{Comparison between the contact rate Kot-type function $\protect%
\beta (t)$ (Black) and the vaccination Kot-type function $v(t)$ (Red) ($%
\protect \varepsilon =0.138$, $\protect \alpha =0.009$, $r=2$ and $\protect%
\phi =\frac{7\protect \pi }{5}$).}
\label{BetaAndVaccination}
\end{figure}

The vaccination switches individuals directly from the susceptible state ($S$%
) to the immune state ($R$). We suppose that infection transmission rate $%
\beta (t)$ is subjected to an environmental periodic forcing and management
imposes a countercycle control strategy through $v\left( t\right) $. The
rationale underpinning a vaccination policy is to ensure that the proportion
of susceptible individuals in the population would stay below a certain
threshold. We display in Fig.~\ref{BetaAndVaccination} a joint
representation of the yearly periodic functions $\beta (t)$ and $v(t)$.
Notice that, for $\phi >0$, the qualitative behavior of $v\left( t\right) $
anticipates the qualitative behavior of $\beta (t).$ In Section III, we will
see how this behavior influences drastically the dynamics, acting as a
controller to suppress chaos. Figure~\ref{BetaAndVaccination} also suggests
that this decisive effect of the phase difference will be achieved with just
a small perturbation, an observable fact reflected by the scale differences
between $\beta (t)$ and $v(t)$ ($v(t)<<\beta (t)$). 

The model is clearly understood when $\beta (t)=$ $\beta _{0}$ and $v(t)=$ $%
0 $. For this case of constant contact rate, and in the absence of a
vaccination strategy, the dynamical system (\ref{SystemVaccination}) has two
equilibrium points:

(i) The Disease-Free Equilibrium (DFE)%
\begin{equation}
\left( S_{0}^{\ast },I_{0}^{\ast },R_{0}^{\ast }\right) =\left( \frac{\sigma 
}{\mu },0,0\right) ,  \label{FixedPoint1}
\end{equation}%
corresponding to a population with no infected individuals;

(ii) The Endemic Equilibrium (EE)%
\begin{equation}
\left( S_{1}^{\ast },I_{1}^{\ast },R_{1}^{\ast }\right) =\left( \frac{\gamma
+\mu }{\beta _{0}},\frac{\mu }{\beta _{0}}\left( \mathcal{R}_{0}-1\right) ,%
\frac{\gamma }{\beta _{0}}\left( \mathcal{R}_{0}-1\right) \right) ,
\label{FixedPoint2}
\end{equation}

corresponding to the case in which there is a significant group of
infectious individuals, where%
\begin{equation}
\mathcal{R}_{0}=\frac{\beta \sigma }{\mu \left( \mu +\gamma \right) }\text{
\  \ }\left( \text{i.e., }\mathcal{R}_{0}=\frac{\beta }{\left( \mu +\gamma
\right) }\text{, for }\sigma =\mu \right)  \label{RoFormula}
\end{equation}%
is a derived basic reproduction number (or basic reproductive number) with
threshold properties. The basic reproduction number of an infection, $R_{0}$%
, is a measure of the potential for the disease to spread in a population.
Heuristically, $R_{0}$\ may be read as the expected number of secondary
infections directly generated by a single infectious individual in a wholly
susceptible population. The number $R_{0}$\ is not a biological constant for
a pathogen as it is also affected by other factors such as the behavior of
the infected population and environmental conditions. It can also be
modified by physical distancing and other public policy or social
interventions. $R_{0}$\ values are usually estimated from mathematical
models, and the estimated values are dependent on the model used and values
of other parameters. Therefore, values given in the literature only make
sense in the given context and it is recommended not to compare directly
values based on different models. The most important uses of $R_{0}$ are
determining if an emerging infectious disease can spread in a population and
determining what proportion of the population should be immunized through
vaccination in order to eradicate a disease. Generally speaking, and
independently from biologically meaningful initial values, this means that:

(a) If $R_{0}<1$, the epidemic cannot maintain itself, since each infected
individual on average infects less than one member of the population. The EE
point $\left( S_{1}^{\ast },I_{1}^{\ast },R_{1}^{\ast }\right) $ is then
unstable, while DFE $\left( S_{0}^{\ast },I_{0}^{\ast },R_{0}^{\ast }\right) 
$ is locally stable and the disease goes extinct;

(b) If $R_{0}>1$, each infected individual infects more than one other
member of the population and a self-sustaining group of infectious
individuals will propagate. In this case, the EE point $\left( S_{1}^{\ast
},I_{1}^{\ast },R_{1}^{\ast }\right) $ is locally stable, while the DFE
point $\left( S_{0}^{\ast },I_{0}^{\ast },R_{0}^{\ast }\right) $ is
unstable. The disease will remain permanently endemic in the population.

Mathematical modeling in epidemiology provides understanding of the
underlying mechanisms that influence the spread of the disease and, in the
process, it suggests control strategies. One of the significant results in
mathematical epidemiology is that the majority of mathematical epidemic
models, including those that have a high degree of heterogeneity, usually
exhibit threshold behavior. In epidemiological terms, this can be stated as
follows: If the average number of secondary infections caused by an average
infectious, the previously called basic reproduction number, is less than
one a disease will die out, while if it exceeds one there will be an
epidemic.\ In the context of differential equation models (or, more
generally, evolution equation models), $R_{0}$ arises as a dimensionless
number of transmission. Throughout this work, we will use our parameters
following previous studies in the literature (see \cite{Zhang2013}, \cite%
{Schwartz1983}, \cite{Kamo2002} and references therein). As a prior notice,
it is important to emphasize that, in the framework of our study, the main
qualitative features of the dynamics are not sensitive to the precise values
of certain parameters. With the time measured in units of years, the values
of the parameters, corresponding to an infectious disease, are listed in the
following table.\pagebreak

\begin{center}
Table 2. List of the parameter values.

$%
\begin{tabular}{|l|}
\hline
$%
\begin{array}{c}
\sigma =\mu =0.01\text{ (year)}^{-1} \\ 
\text{ }%
\end{array}%
\begin{array}{l}
\Longrightarrow \frac{1}{\mu }=100\text{ years} \\ 
\text{(mean lifetime of the host)}%
\end{array}%
$ \\ \hline
$%
\begin{array}{c}
\gamma =50\text{ (year)}^{-1} \\ 
\text{ }%
\end{array}%
\begin{array}{l}
\Longrightarrow \frac{1}{\gamma }=0.02\text{ year}\approx \text{1 week} \\ 
\text{ (mean infectious period)}%
\end{array}%
$ \\ \hline
$%
\begin{array}{c}
\beta _{0}=1505\text{ (year)}^{-1} \\ 
\text{ }%
\end{array}%
\begin{array}{l}
\Longrightarrow \frac{1505}{365}\approx 4\text{ (day)}^{-1} \\ 
\text{ (mean number of contacts per day)}%
\end{array}%
$ \\ \hline
$%
\begin{array}{c}
R_{0}\approx \frac{\beta _{0}}{\left( \mu +\gamma \right) } \\ 
\text{ }%
\end{array}%
\begin{array}{l}
\approx 30 \\ 
\text{ (average number of secondary infections)}%
\end{array}%
$ \\ \hline
\end{tabular}%
$
\end{center}

Within proper meaningful contexts, the degree of seasonality $\varepsilon $,
of the main force $\beta \left( t\right) $, as well as parameters $v_{0}$, $%
\alpha $ and $\phi $ of the vaccination rate (the added small perturbation
that will act as a control component), will be taken as control parameters.

At this moment, it is also critical to stress that, throughout our study, a
close attention will be focused, not on the phase control method \textit{per
se}, but on its biological significance and on its consequences in an
epidemiological context.

\section{Existence of chaos: epidemics with a seasonal contact rate and
constant vaccination}

\begin{figure*}[tbp]
\centerline{\includegraphics[width=0.45\textwidth]{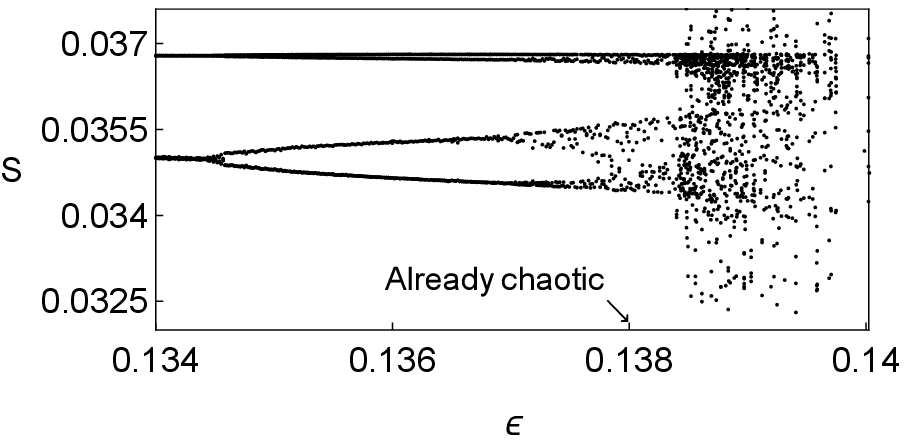} 
                  \hspace{0.15cm}
                  \includegraphics[width=0.45\textwidth]{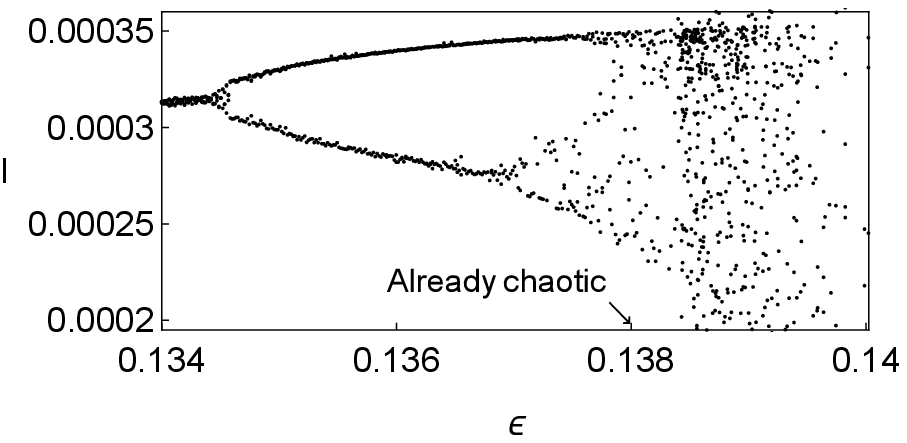}}%
\centerline{\includegraphics[width=0.45\textwidth]{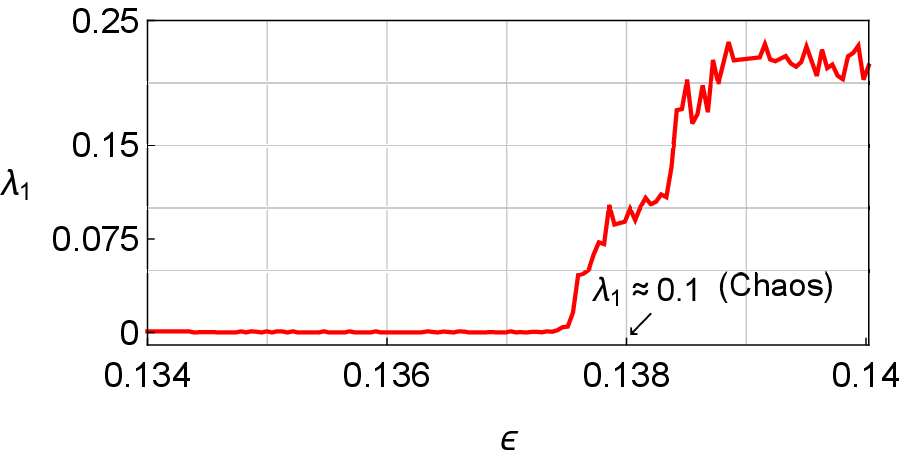}}%
\centerline{\includegraphics[width=0.45\textwidth]{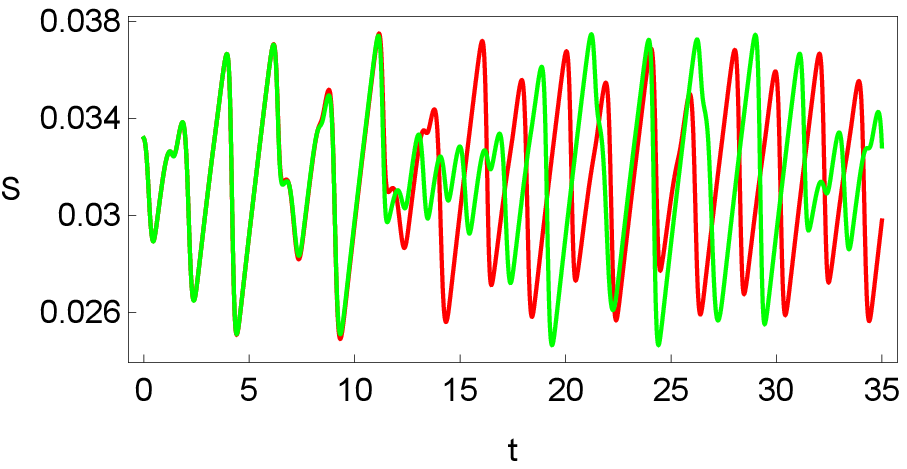} 
                  \hspace{0.15cm}
                  \includegraphics[width=0.45\textwidth]{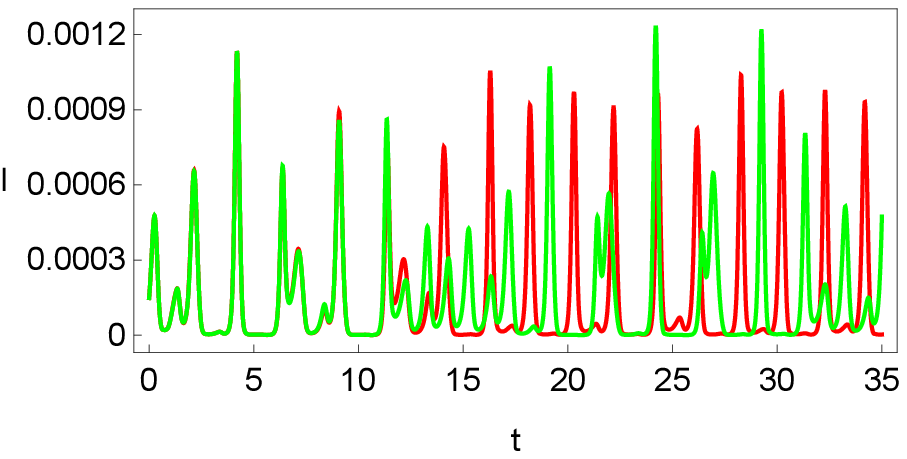}}
\caption{Global noticeable features of the dynamics of the SIR system (%
\protect \ref{SystemVaccination}) considering a constant vaccination
strategy, $v\left( t\right) =v_{0}=0.071$ and $\protect \alpha =0$.
Increasing values of the degree of seasonality of the contact rate, $\protect%
\varepsilon $, bring the dynamics to a chaotic regime. Upper panel:\textbf{\ 
}bifurcation diagrams of the dynamical variables $S$ and $I$, taking $%
\protect \varepsilon $ as control parameter ($0.134\leq \protect \varepsilon %
\leq 0.14$). Middle panel: variation of the largest Lyapunov exponent, $%
\protect \lambda _{1}$, with $\protect \varepsilon $ ($0.134\leq \protect%
\varepsilon \leq 0.14$). Lower panel: corresponding time series of the
dynamical variables $S$\ and $I$, obtained using two different (but close)
initial conditions displayed in \textit{red} and \textit{green} ($\protect%
\varepsilon =0.138$, value for which the system is chaotic). }
\label{Bifs_Lyap}
\end{figure*}

The complexity of the dynamics increases considerably with the introduction
of the seasonal contact rate. In this section, we will examine the long-term
behavior of the three-dimensional chaotic attractors arising in the forced
SIR system (\ref{SystemVaccination}), with a seasonal transmission rate
component and under a constant vaccination strategy. The time series,
displayed in Fig.~\ref{Bifs_Lyap} (Lower panel), govern the population
dynamics of the susceptible and infected. They were obtained using two
different, but close, initial conditions.

In this paragraph, the Lyapunov exponents of the SIR model (\ref%
{SystemVaccination}) with the seasonal component $\beta \left( t\right) $
and $v\left( t\right) =v_{0}$, receive our attention as a framework to
diagnose chaos in the system. A discussion about the Lyapunov exponents as a
quantitative measure of the rate of separation of infinitesimally close
trajectories, as well as a computation method, can be found in \cite%
{ParkerChua1989}. The characteristic Lyapunov exponents measure the typical
rate of the exponential divergence of nearby trajectories in phase space,
i.e., they give us information on the rate of growth of a very small error
on the initial state of the system. 

A positive Lyapunov exponent is commonly taken as an indicator of chaotic
behavior. In Fig.~\ref{Bifs_Lyap}, taking $\varepsilon $ as a control
parameter, we present the variation of the largest Lyapunov exponent $%
\lambda _{1}$ (Middle panel) and two\ bifurcation diagrams (Upper panel)
regarding the dynamical variables $S$ and $I$. In agreement with the
represented bifurcation diagrams, the largest Lyapunov exponent, $\lambda
_{1}$, is positive in the chaotic regime and $\lambda _{1}\approx 0$ in the
periodic windows (precisely in agreement with the theory for 3D systems).
When this Lyapunov exponent is positive, it indicates the region where the
system is chaotic. Increasing the degree of seasonality $\varepsilon $
results in the complexity of the dynamics at higher values. This way, the
actual existence of chaos is numerically recognized, with the seasonal
transmission function $\beta \left( t\right) $ as the force driving the
system to the chaotic state. In particular, for the choice of parameters $%
\beta _{0}=1505$, $\varepsilon =0.138$ of the contact rate $\beta \left(
t\right) $ and $v\left( t\right) =v_{0}=0.071$, the system is chaotic and
the largest Lyapunov exponent is $\lambda _{1}\approx 0.1$. From now on, we
keep the parameter values ($\beta _{0}=1505$, $\varepsilon =0.138$ and $%
v_{0}=0.071$) all throughout our study.

\section{Suppressing chaos: the phase effect on a seasonal vaccination
strategy}

As we have just stated, intense seasonality of the periodically varying
contact rate induces chaotic dynamics in the epidemic system under a
conventional constant vaccination strategy. However, in the context of
epidemiology, chaos is often regarded as an undesirable phenomenon
associated with the erratic permanence of infectious diseases. 
\begin{figure*}[tbp]
\centerline{\includegraphics[width=0.355\textwidth]{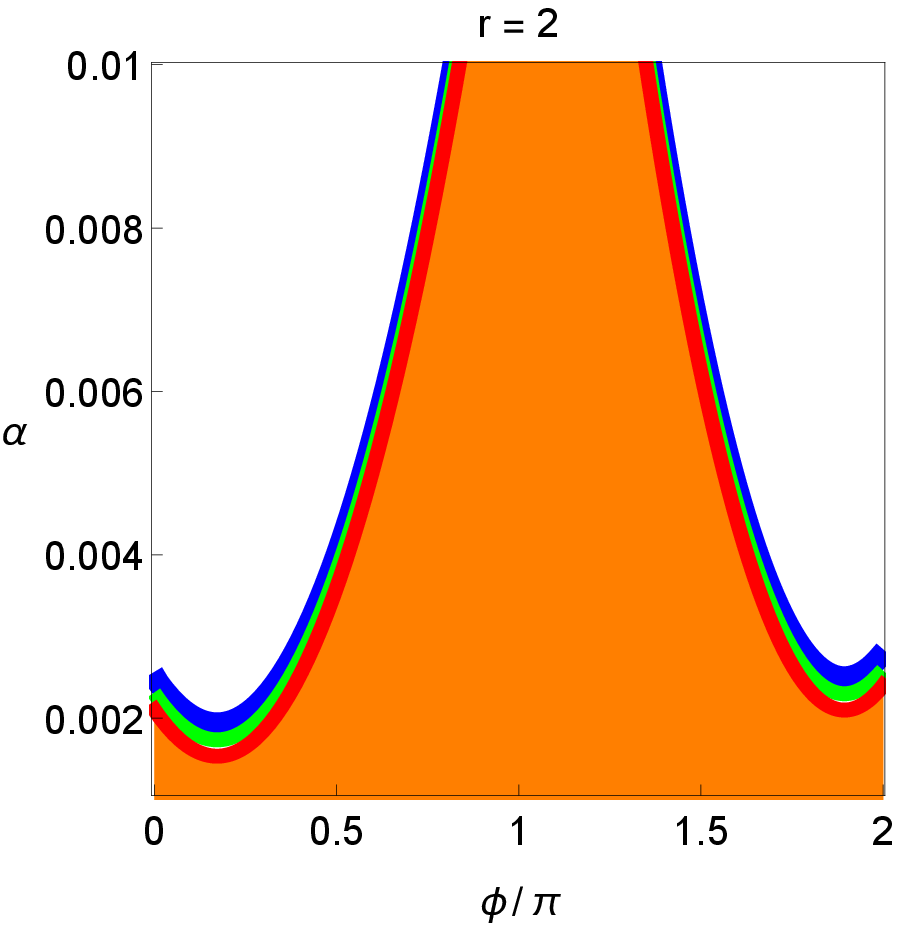} 
                  \hspace{0.15cm}
                  \includegraphics[width=0.5\textwidth]{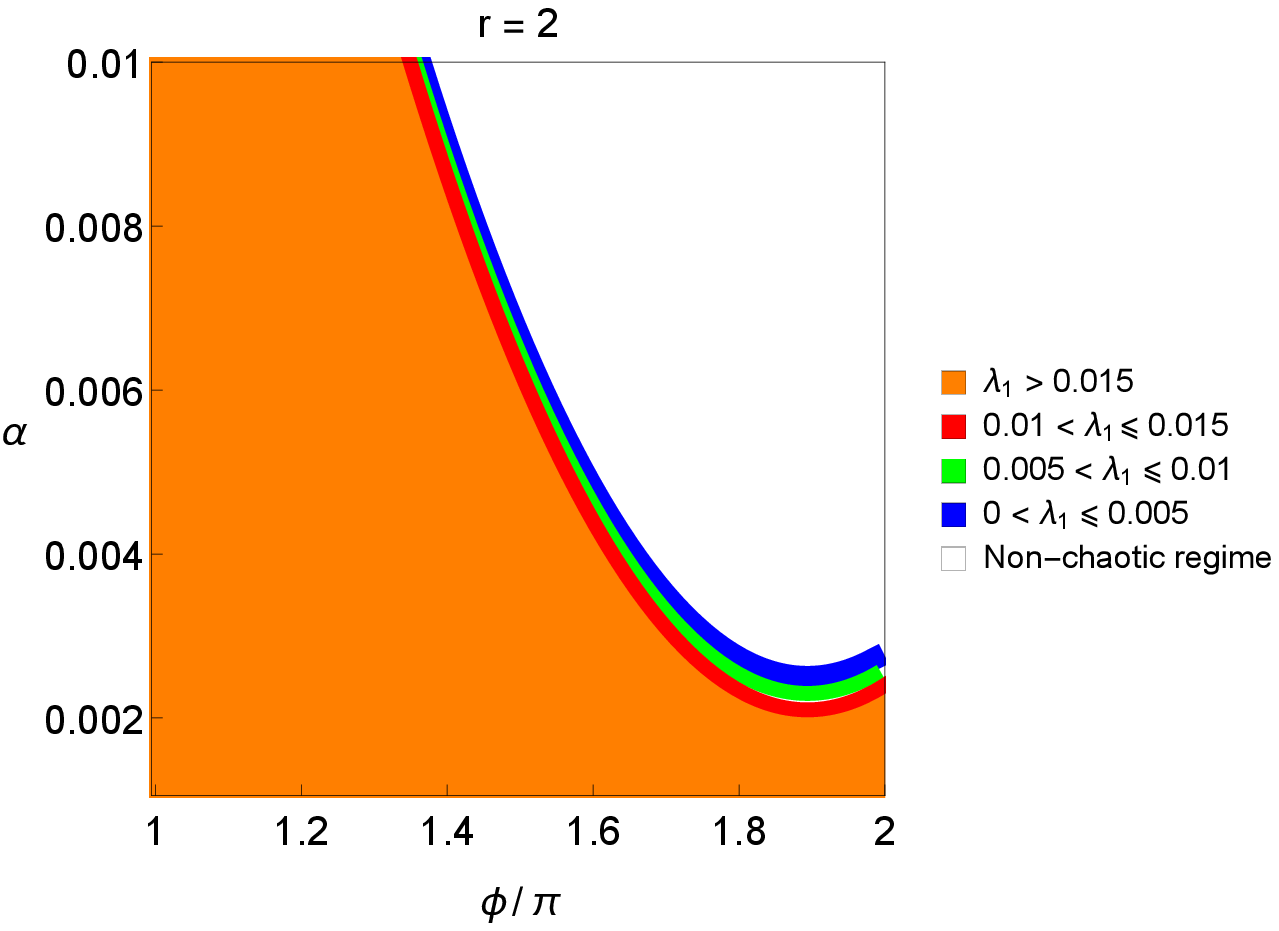}}
\caption{Density plots for the largest Lyapunov exponent corresponding to $%
r=2$, in the $\left( \protect \phi ,\protect \alpha \right) $-parameter
region. Left: $0\leq \protect \phi \leq 2\protect \pi $ and $0.001\leq \protect%
\alpha \leq 0.01$; Right: \textit{zoom in} of the previous density plot
considering $\protect \pi \leq \protect \phi \leq 2\protect \pi $ and $%
0.001\leq \protect \alpha \leq 0.01$. We distinguish two dynamical regimes:
the chaotic regime, corresponding to $\protect \lambda _{1}>0$ (colored
region) and the regular regime, represented by the white region.}
\label{Density_Plots}
\end{figure*}
As a consequence, the aim of the present section is to provide a
comprehensive study of the forced epidemic system in terms of the
implementation of the phase control, as a control technique to suppress the
chaotic behavior.

Inspired by the importance of vaccination for the elimination of infectious
diseases (please see \cite{Pasour2012}), this procedure is introduced as a
biologically meaningful weak perturbation on a seasonal vaccination
strategy. More specifically, with the parameters $\beta _{0}=1505$, $%
\varepsilon =0.138$ and $v_{0}=0.071$ already tailored in such a way that
the asymptotic state of system (\ref{SystemVaccination}), under a constant
vaccination strategy ($v\left( t\right) =v_{0}$), is chaotic, our aim here
is to analyze the effects of the phase $\phi $ in the chaotic regime, when
an additive perturbation to $v_{0}$, given by (\ref{Additive_Perturbation}),
is included. For this purpose, the amplitude of the perturbation $\alpha $
is assumed to be very small, i.e., $\alpha <<1$.

\begin{figure*}[tbp]
\centerline{\includegraphics[width=0.45\textwidth]{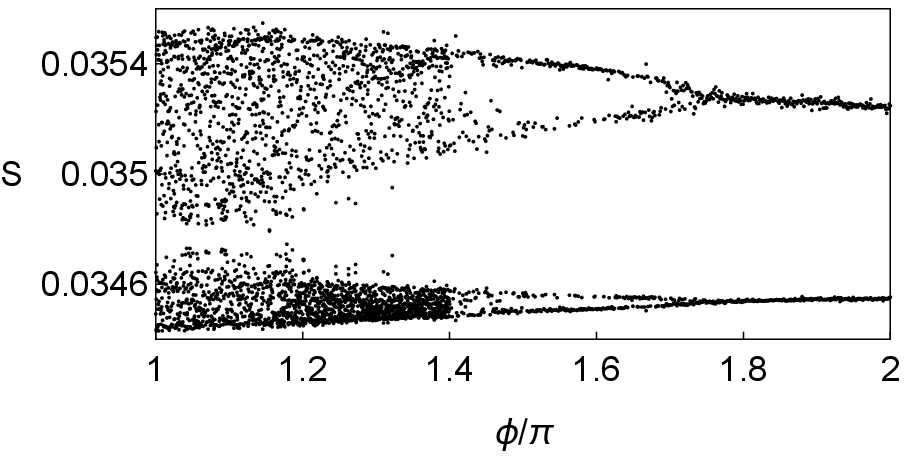} 
                  \hspace{0.15cm}
                  \includegraphics[width=0.45\textwidth]{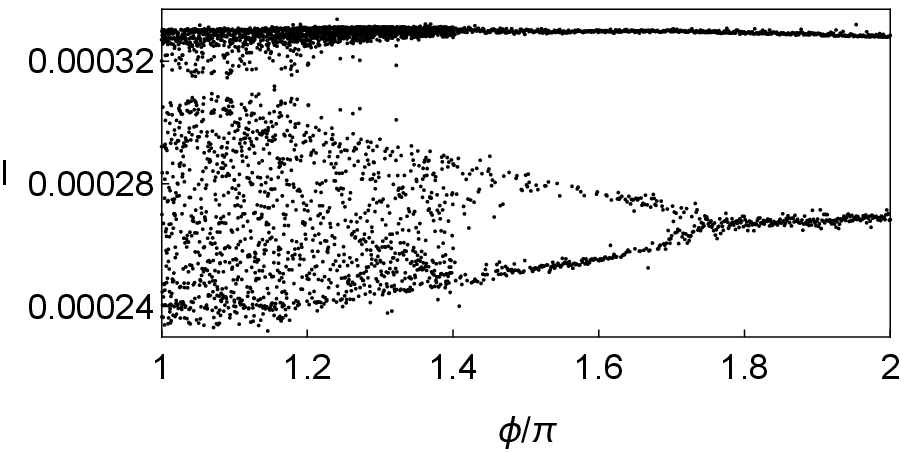}}%
\centerline{\includegraphics[width=0.45\textwidth]{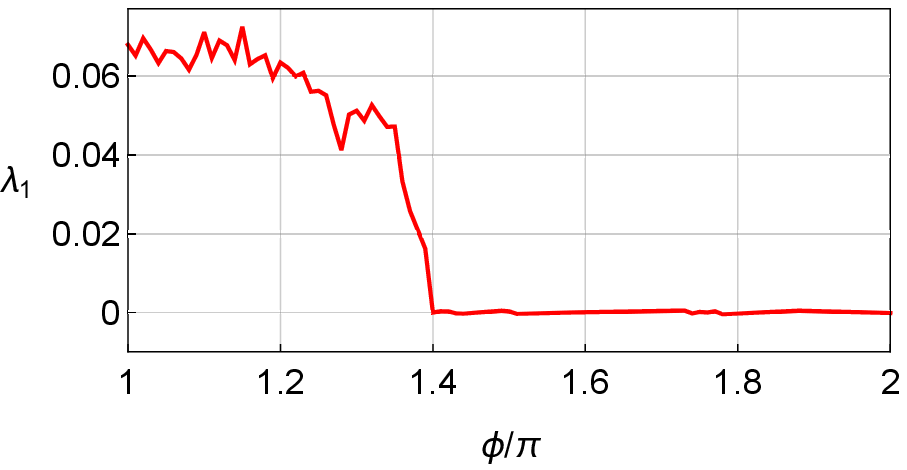}}%
\centerline{\includegraphics[width=0.45\textwidth]{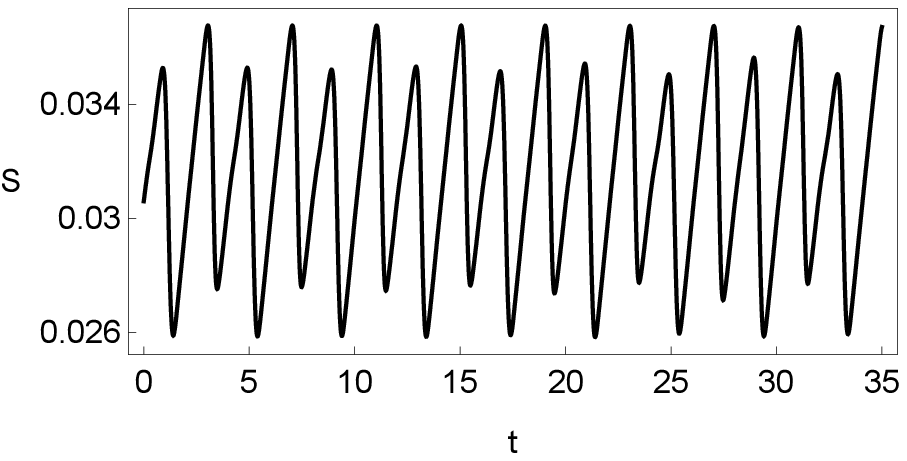} 
                  \hspace{0.15cm}
                  \includegraphics[width=0.45\textwidth]{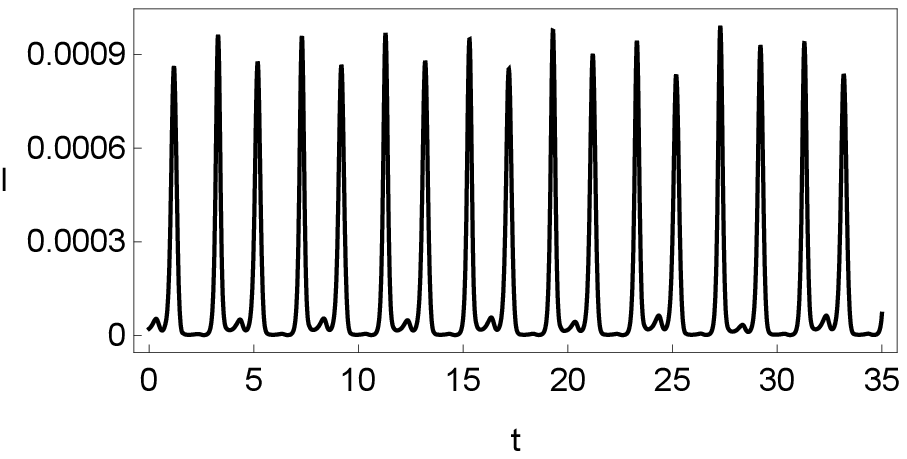}}
\caption{Upper panel: bifurcation diagrams for $S$ and $I$ when $\protect%
\alpha =0.009$, taking $\protect \phi $ as control parameter ($\protect \pi %
\leq \protect \phi \leq 2\protect \pi $); Middle panel: variation of the
largest Lyapunov exponent, with $\protect \alpha =0.009$ and $\protect \pi %
\leq \protect \phi \leq 2\protect \pi $; Lower panel: time series for the
dynamical variables $S$ and $I$, for $\protect \alpha =0.009$ and $\protect%
\phi =\frac{7}{5}\protect \pi $ (a couple of parameter values marks the
beginning of a transition from chaotic to regular behavior).}
\label{Phase_Control}
\end{figure*}
In order to visualize the effect of the phase $\phi $, combined with the
variation of $\alpha $, we have displayed in Fig.~\ref{Density_Plots}
density plots. The color scale corresponds to the computation of the largest
Lyapunov exponent over the parameter region characterized by $0\leq \phi
\leq 2\pi $ and $0.001\leq \alpha \leq 0.01$, fixing $r=2$. In Fig. \ref%
{Density_Plots}, we distinguish two dynamical regimes. The chaotic regime,
corresponding to $\lambda _{1}>0$, represented by the colored region (Blue: $%
0<\lambda _{1}\leq 0.005$; Green: $0.005<\lambda _{1}\leq 0.01$; Re\.{d}: $%
0.01<\lambda _{1}\leq 0.015$; Orange: $\lambda _{1}>0.015$) and the regular
regime, represented by the white region. Taking $r=2$, we are able to obtain
a clear transition from chaotic to regular behavior increasing $\phi $. In a
biological context, this choice of $r$ guarantees a frequency of the
vaccination within one year. With these density plots, we can better
appreciate the structure of the chaotic region (colored) and periodic region
(white). In our numerical exploration, the choice of such small values for
the perturbation parameter $\alpha $ allows us to demonstrate the
effectiveness of the phase control. The density plot (right-hand side of
Fig.~\ref{Density_Plots}), is precisely a zoom in of the previous density
plot (left-hand side of Fig.~\ref{Density_Plots}), which clearly exhibits
the suppression of chaos, increasing the values of $\phi $, when $\pi \leq
\phi \leq 2\pi $. In particular, the striking characteristics of the phase
control - the key role of the phase $\phi $ in selecting the final state of
the system successfully, combined with tiny values of $\alpha $ to suppress
chaos - are present in this density plot (right-hand side of Fig.~\ref%
{Density_Plots}). In order to continue illustrating this primary role of the
phase $\phi $ we show, in Fig.~\ref{Phase_Control} (upper panel),
bifurcation diagrams for $S$ and $I$, fixing $\alpha =0.009$ and taking $%
\phi $ as control parameter, with $\pi \leq \phi \leq 2\pi $ and the
corresponding variation of the largest Lyapunov exponent (Middle panel). In
fact, we have found a wide range of phase values producing regular motion,
in which the system leaves the chaotic region to periodic motion via inverse
period doubling. A remarkable feature is that for values of $\alpha $, such
that $\alpha \succsim 0.003$, chaos is entirely wiped out in this phase
region after having achieved a certain threshold of the phase $\phi $
(density plot of Fig.~\ref{Density_Plots}, right-hand side). In particular,
for $\alpha =0.009 $, the periodic state is obtained when $\phi \approx 
\frac{7}{5}\pi $. This observation indicates that the phase $\phi $ is a
sensitive parameter for the system bifurcation, i.e., the distribution of
regular and chaotic regions strongly depends on the phase difference
(please, visit again Fig. \ref{BetaAndVaccination}, and see the
representation of $v\left( t\right) $ precisely for $\phi \approx \frac{7}{5}%
\pi $. It is extremely interesting that the weak additive force makes the
vaccination strategy seasonal and its phase term (with $\phi >0$) ensures
that the vaccination $v\left( t\right) $ anticipates the seasonality of the
transmission rate $\beta \left( t\right) $. Therefore, the phase control
method for suppressing chaos takes on a clear meaning, where the phase
difference between the two central driven forces: (i) a transmission rate
and (ii) a vaccination strategy, stand out to be a key perturbation
parameter with immediate and beneficial biological consequences - the
control of a given infectious disease.

\section{Concluding remarks}

For decades, mathematical models have been proposed to evaluate the spread
and control of infectious diseases. Recently, for the last entire year, the
world has been experiencing the intensification of the research in epidemic
dynamics generated by the dissemination of COVID-19 (\cite{Carlson2020},\cite%
{Wu2020},\cite{Lin2020},\cite{Guan2020} and \cite{Zhu2020}). The year 2020
has seen significant advances taking place, to build the infrastructures to
keep up with this coronavirus. SARS-CoV-2, the virus which causes COVID-19,
is constantly evolving and mutating, as do all similar virus. The pressure
on the virus to evolve is increased by the fact that so many millions of
people have now been infected. Such identified changes since it emerged in
2019, in a RNA virus that exists as a cloud of genetic variants known as
quasispecies, are completely to be expected to occur and have been useful in
understanding the worldwide spread as well as the transmission patterns. The
majority of the mutations will not be significant or cause for concern, but
some may give the virus an evolutionary advantage which may lead to higher
transmission. More specifically, a new genetic variant of the virus has
emerged and is spreading in many parts of the UK and across the world. This
efficient transmission among people is usually associated to a modeling
process with a higher basic reproduction number, $R_{0}$.

It is precisely in this context of controlling infectious diseases with $%
R_{0}>>1$, that we have successfully applied the phase control method of
suppressing chaos to the continuous periodically driven SIR epidemic model
with a seasonal transmission rate and under a conventional constant
vaccination strategy. A close attention has been devoted, not to the phase
control procedure \textit{per se}, but to its biological significance and to
its consequences within the framework of epidemiology. We have provided
detailed/revealing insights about the role played by the phase difference of
the two periodically driven forces - the seasonal transmission function $%
\beta \left( t\right) $ and a vaccination component $v\left( t\right) $.
Having started with the analysis of the existence of chaos, the transmission
rate function $\beta \left( t\right) $ has been identified as the main force
driving the system to the chaotic regime under a classical constant
vaccination strategy. With the computation of the largest Lyapunov exponent,
values of the parameters of $\beta \left( t\right) $ have been tailored
within the chaotic region.

Motivated by the fact that, in an epidemiological context, chaotic behavior
is often associated with the erratic permanence of infectious diseases and a
vaccination strategy is associated with their efficient elimination. Thus,
an idea emerged - to introduce the control procedure as a biologically
meaningful weak perturbation on a seasonal vaccination function $v\left(
t\right) $. Given the importance of controlling the chaotic behavior, i.e.,
the necessity of having predictable densities for the epidemic populations,
we have applied the mentioned periodic control signal $v\left( t\right) ,$
including the phase difference with respect to the periodic forcing of the
initial system, which has acted as an effective control strategy. More
precisely, the chaotic epidemic outbreaks, that appeared as a result of the
seasonal variations in the contact rate, have been suppressed by the used
vaccination control scheme. Indeed, the crucial role of the phase term, in
the seasonal component of $v\left( t\right) $, was evidenced by using
numerical simulations that allowed us to clearly identify dynamical
transitions from a chaotic regime to a regular behavior, which is
biologically associated with the control of a given infectious disease.

This study provides another illustration of how an integrated approach,
involving numerical evidences and theoretical reasoning, within the theory
of dynamical systems, can contribute to our understanding of important
biological models and provide a trustworthy explanation of complex phenomena
witnessed in biological systems.

Above all, the recent appearance of COVID-19 is a reminder that there is
still so much to learn about a pandemic dynamics. The pace of the research
effort in the past year has been extraordinary. However, there is no room
for complacency. We have to be humble and mostly be prepared to adapt and
respond to new and continued changes.

\textbf{Acknowledgments}

This research has been financially supported by the Portuguese Foundation
for Science and Technology (FCT) under Projects No.~UIDB/04106/2020 and
No.~UIDP/04106/2020 (CIDMA) (CJ) and No.~UID/MAT/04459/2013 (JD and NM); and
by the Spanish State Research Agency (AEI) and the European Regional
Development Fund (ERDF, EU) under Projects No.~FIS2016-76883-P and
No.~PID2019-105554GB-I00 (JMS and MAFS).

\end{document}